\documentclass[11pt]{article}
\usepackage{amsmath}
\usepackage{amsfonts}
\usepackage{amssymb}
\usepackage{enumerate}
\usepackage{tabularx,booktabs}
\usepackage{float}
\usepackage{graphicx}
\usepackage{url}

\newtheorem{theorem}{Theorem}
\newtheorem{conjecture}[theorem]{Conjecture}

\newcommand{\QED}{\hfill$\square$}

\title {
    \bf {Counterexamples to conjectures on graph distance measures based on topological indexes}
}

\author
{
{\large \sc Aleksandar Ili\' c } \\
{\em \normalsize Facebook Inc, Menlo Park, CA, USA} \\
{\normalsize e-mail: { \tt aleksandari@gmail.com }}
\and
{\large \sc Milovan Ili\' c } \\
{\em \normalsize Metropolitan University, Belgrade, Serbia} \\
{\normalsize e-mail: { \tt ilicmilovan@gmail.com }}
}

\begin{document}

\maketitle

\begin{abstract}
In this paper we disprove three conjectures from [M. Dehmer, F. Emmert-Streib, Y. Shi, \emph{Interrelations of graph distance measures based on topological indices}, PLoS ONE 9 (2014) e94985] on graph distance measures based on topological indices by providing explicit classes of trees that do not satisfy proposed inequalities. The constructions are based on the families of trees that have the same Wiener index, graph energy or Randi\' c index - but different degree sequences.
\end{abstract}

\medskip

\medskip

\section{Introduction}

The structural graph similarity or distance of graphs has been attracting attention of researchers from many different fields, such as mathematics, social network analysis, biology or chemistry. The two main concepts which have been explored are exact (based on isomorphic relations and computationally demanding, see \cite{DeMe07, DeEm07} and references therein) and inexact graph matching (one example here is graph edit distance, see \cite{Bu00}). In this paper we are dealing with the later concept, and study interrelations of graph distance measures by means of inequalities. Graph distance measure is a mapping $d:\mathcal{G} \times \mathcal{G} \rightarrow R^+$, where $\mathcal{G}$ is any set of graphs \cite{EmDeSh16}.

Given two graphs with $n$ vertices, $G_1$ and $G_2$, the edit distance between $G_1$ and $G_2$, denoted by $GED(G_1,G_2)$, is the minimum cost caused by the number of edge additions and/or deletions that are needed transform $G_1$ into $G_2$. $GED$ has been introduced by Bunke \cite{Bu83}.

A topological index \cite{ToCo02} is a type of a molecular descriptor that is calculated based on the molecular graph of a chemical compound, and usually represents numerical graph invariant $I(G)$. 

Some novel class of graph distance measures based on topological indices have been introduced by Dehmer et al. \cite{DeEmSh14}. The same authors continued research from \cite{DeEmSh15} and proved
various inequalities and studied comparative graph measures based on the well-known Wiener index, graph energy and Randi\' c index. For more results on graph entropy measures see \cite{CaDeSh14,DeMo11}.

\section{Preliminaries}

We focus on the graph distance measure introduced in \cite{DeEmSh14}:
$$
d_I(G, H) = d(I(G), I(H)) = 1 - e^{-\left(\frac{I(G) - I(H)}{\sigma}\right)^2},
$$
where $\sigma$ is an arbitrary real number and $I(G)$ and $I(H)$ are certain graph invariants. Note that $d_I(G, H)$ is actually is a distance measure for real numbers, and for more mathematical properties see \cite{DeEmSh15}.

\bigskip

Let $G = (V, E)$ be a connected graph on $n$ vertices. The distance between the vertices $u$ and $v$ of $G$ is denoted by $d(u, v)$. The Wiener index is one of the oldest distance-based topological indexes \cite{DoEn01,IlIl13}, and defined as the sum of all distances between any two vertices of a graph:
$$
W(G) = \sum_{u, v \in V} d(u, v).
$$
Randi\' c index \cite{LiSh08,RoTo15} is suitable for measuring the extent of branching of chemical graphs, and is defined as
$$
R(G) = \sum_{uv \in E} \frac{1}{\sqrt{deg(u) deg(v)}},
$$
where $deg(v)$ is a degree of the vertex $v$.

Denote by $\lambda_1, \lambda_2, \ldots, \lambda_n$ the eigenvalues of the adjacency matrix of $G$. The energy of a graph is introduced by Gutman \cite{LiSh12} equals to
$$
E(G) = \sum_{i = 1}^n |\lambda_i|,
$$
while the distance measure based on Shannon's entropy \cite{IlDe15} is defined as
$$
Ig(G) = \log E(G) - \frac{1}{E(G)} \sum_{i = 1}^n |\lambda_i| \log |\lambda_i|.
$$

By setting $f(v_i) = deg(v_i)^k = deg_i^k$ in the Shannon entropy formula, we can also obtain the new entropy based on the degree powers \cite{CaDeSh14,Il16}, denoted by $If_k(G)$:
$$
If_k(G) = \log \left (\sum_{i = 1}^n deg_i^k \right) - \frac{1}{\sum_{i = 1}^n deg_i^k} \sum_{i = 1}^n deg_i^k \log deg_i^k.
$$

The following three conjectures are proposed in the same paper of Dehmer et al.

\begin{conjecture}
Let $T$ and $T'$ be any two trees with $n$ vertices. Then, it holds
$$
d_W (T, T') \geq d_R(T, T').
$$
\end{conjecture}

\begin{conjecture}
Let $T$ and $T'$ be any two trees with $n$ vertices. Then, it holds
$$
d_E (T, T') \geq d_{Ig} (T, T').
$$
\end{conjecture}

\begin{conjecture}
Let $T$ and $T'$ be any two trees with $n$ vertices. Then, it holds
$$
d_R  (T, T') \geq d_{If_1} (T, T')
$$
\end{conjecture}

In this note, we are going to construct a family of pairs of trees that do not satisfy Conjecture 1 and 3, and then disprove Conjecture 2 by providing specific examples of equienergetic trees from the literature. The authors from \cite{DeEmSh14} verified the conjectures on all trees with small number of vertices and these counterexamples have more than 12 vertices. In Section 4, we are going to refine some results on the graph edit distance and Shannon entropy from \cite{DeEmSh14} and \cite{DeEmSh15}.

\section{Main results}

By definition, the reverse inequality $d_W (T, T') < d_R(T, T')$ is equivalent to
$$
1 - e^{-\left(\frac{W(T) - W(T')}{\sigma}\right)^2} < 1 - e^{-\left(\frac{R(T) - R(T')}{\sigma}\right)^2}.
$$
Using the fact that $e^x$ is a strictly increasing function, it is further equivalent to
$$
|W(T) - W(T')| < |R(T) - R(T')|.
$$
We are going to construct a family of tree pairs $(T, T')$ on $n$ vertices that have the same Wiener index and different Randi\' c index, and thus left side of the above inequality will be 0. More examples of non-isomorphic trees having equal Wiener index are introduced by Rada in~\cite{Ra05}.

\begin{figure}[h]
\medskip
\centering
\includegraphics[scale=0.2]{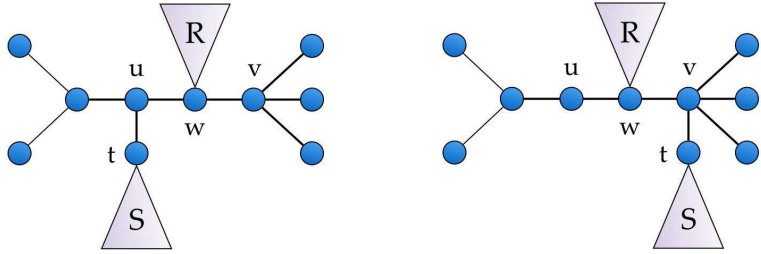}
\caption{Trees $T$ and $T'$ with equal Wiener index.}
\medskip
\end{figure}

In Figure 1, we can attach arbitrary trees $S$ and $R$ to the vertices $t$ and $w$, and it is an easy exercise to verify that these trees are non-isomorphic and have the same Wiener index (the sum of all distances from the vertex $t$ to the vertices not in $S$ and $R$ is constant in both trees and equals 26). 

Let $deg(t) = x$ and $deg(w) = y$. Randi\' c index is a local topological index, and the difference between $R(T)$ and $R(T')$ depends only on the edge weights around vertices $u$, $w$, $v$ and $t$. By direct computation, it follows 
\begin{eqnarray*}
R(T) - R(T') &=& \left ( \frac{1}{\sqrt{3 \cdot 3}} + \frac{1}{\sqrt{3 \cdot x}} + \frac{1}{\sqrt{3 \cdot y}}
+ \frac{3}{\sqrt{1 \cdot 4}} + \frac{1}{\sqrt{4 \cdot y}} \right) \\
&& - \left( \frac{1}{\sqrt{3 \cdot 2}} + \frac{1}{\sqrt{2 \cdot y}} + \frac{1}{\sqrt{5 \cdot y}} + \frac{3}{\sqrt{1 \cdot 5}} + \frac{1}{\sqrt{5 \cdot x}} \right)  
\end{eqnarray*}

As $x$ and $y$ can be arbitrary integer numbers, the above difference is almost never equal to zero - and this completes the counterexample construction for Conjecture~1.

For Conjecture 3, we need to construct a family of pairs of trees that satisfy
$$
|R(T) - R(T')| < |If_k(T) - If_k(T)|,
$$
for any $k \geq 1$. Let $C(x, y, z, t, T)$ be a tree composed of caterpillar having $x$, $y$, $z$ and $t$ leaves attached to any tree $T$. Then by definition it holds
\begin{eqnarray*}
R(C(x, y, z, t, T)) &=& \frac{x-1}{\sqrt{x}} + \frac{y-2}{\sqrt{y}} + \frac{z-2}{\sqrt{z}} + \frac{t-2}{\sqrt{t}} \\
&& + \frac{1}{\sqrt{xy}} + \frac{1}{\sqrt{yz}} + \frac{1}{\sqrt{zt}} + R'(T),
\end{eqnarray*}
where $R'(T)$ is the remainder of the summation for the Randi\' c index of the subtree $T$.

We are interested in the pairs of quadruples $(x, y, z, t)$ and $(x', y', z', t')$ such that 
$$
t = t' \qquad \mbox{and} \qquad R(C(x, y, z, t, T)) = R(C(x', y', z', t', T)).
$$

The main trick for the exhaustive search is to fix $t = 4$ and search among perfect squares, as that will make the calculation easier. After running a simple code for iterating through all triples less than 100, we get multiple examples of pairs of trees with equal Randi\' c index (see Figure 2):
$$
R(C(9, 4, 9, 4, T)) = R(C(4, 16, 4, 4, T))
$$
or
$$
R(C(36, 36, 4, 4, T)) = R(C(64, 9, 9, 4, T)).
$$

\begin{figure}[h]
\medskip
\centering
\includegraphics[scale=0.2]{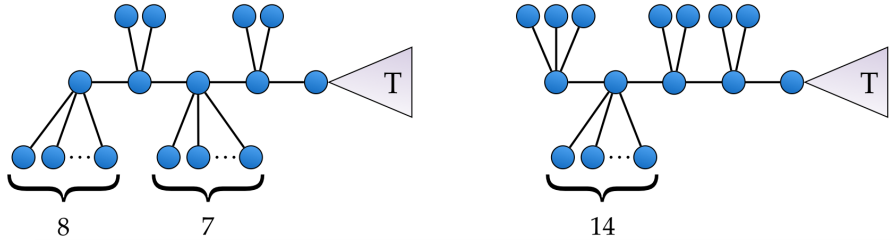}
\caption{Trees $T$ and $T'$ with equal Randi\' c index.}
\medskip
\end{figure}

For these graph pairs, we obviously have $If_1(T) \neq If_1(T')$, given that the degree sequences are different. By definition we get
$$
If_1(G) = \log 2m - \frac{1}{2m} \sum_{i = 1}^n deg_i \log deg_i,
$$
and $If_1(T) \neq If_1(T')$ is equivalent to 
$$
9 \log 9 + 4 \log 4 + 9 \log 9 + 4 \log 4 \neq 4 \log 4 + 16 \log 16 + 4 \log 4 + 4 \log 4.
$$

The above approadh can be easily generalized to $If_k(T) \neq If_k(T')$ for an arbitrary $k \geq 1$.

\bigskip

The same idea can be applied to Conjecture 2, which is equivalent to
$$
|E(T) - E(T')| < |Ig(T) - Ig(T')|.
$$

For trees with the same spectra it obviously holds that $E(T) = E(T')$ and $Ig(T) = Ig(T')$. Here we can construct non-cospectral trees with the same energy and non-zero difference of their $Ig$ invariants.

There are many examples of non-cospectral equienergetic trees in the literature (see \cite{BrSt04,MiFu09} and references therein), and by direct verification for some of them we can conclude that $Ig(T) \neq Ig(T')$. For the first pair of trees defined in \cite{MiFu09}, it can be directly computed that 
$Ig(T_1) = 2.67176585647... < Ig(T_3) = 2.72814227465...$ (by ignoring zero eigenvalues). Note that until now no general method for construction of such trees is known. This refutes Conjecture~2.

\section{Further results}

Let $p=(p_1, p_2, \ldots, p_n)$ be a probability vector, namely $0 \leq p_i \leq 1$ and $\sum_{i = 1}^n p_i = 1$. The Shannon's entropy of $p$ has been defined by
$$
I(p) = - \sum_{i = 1}^n p_i \log p_i.
$$
We denote by $d_{Ip}$ the graph distance measure of $I(p)$. The authors in \cite{DeEmSh14} proved the following result:

\begin{theorem}
Let $G$ and $H$ be two graphs with the same vertex set. Denote by $p = (p_1, p_2, \ldots, p_n)$ and $p' = (p_1', p_2', \ldots, p_n')$ be the probability vectors of $G$ and $H$, respectively. If $p_i \leq p_i'$ for each $1 \leq i \leq n$, then we infer
$$
d_{Ip}(G, H) < 1 - e^{-A^2/\sigma^2},
$$
where
$$
A = \sum_{i = 1}^n \left ( p_i' \log \left(1 + \frac{1}{p_i'} \right) + \log(p_i' + 1) \right).
$$
\end{theorem}

Since the sum of the coordinates of the probability vectors $p$ and $p'$ is equal to 1 and $p_i \leq p_i'$ we easily conclude that $p = p'$ and therefore $d_{Ip} (G, H) = 0$. This finally means that the above result needs to be reformulated by dropping the condition for probability vectors.

\bigskip

The authors from \cite{DeEmSh15} focused on graphs whose distance can be obtained by performing only one graph edit operation, i.e. $GED(G, H) = 1$. This is used to determine which measure is the powerful and most useful one. They presented such results on Randi\' c index and graph energy, but did not obtain similar bounds for the Wiener index. 

If an edge $e = vu$ is removed, the distance between $u$ and $v$ increases and no
other distances decrease. Favaron,
Kouider and Mah\' eo \cite{FaKo89} proved the following result:

\begin{theorem}
Let $G$ be a graph of order $n$ and $e$ a cyclic edge of $G$. Then
$$
\mu(G-e) - \mu(G) \leq \frac{\sqrt{2} - 1}{3} n + O(1),
$$
where $\mu(G)$ denotes the average distance, defined as $\mu(G) = \frac{W(G)}{{n \choose 2}}$.
\end{theorem}

From the previous theorem, we can easily conclude the following relation
$$
1 \leq W(G-e) - W(G) \leq \frac{\sqrt{2} - 1}{6} n^3 + O(n^2). 
$$
 
We complete the section about graph distance measures based on the Wiener index:

\begin{theorem}
Let $G$ and $H$ be two connected graphs with $n$ vertices. If $GED(G, H) = 1$, then 
$$
d_W(G, H) \leq 1 - e^{-\left( \frac{\sqrt{2} - 1}{6} n^3 + O(n^2) \right)^2 / \sigma^2}.
$$ 
\end{theorem}

\bigskip {\bf Acknowledgment. } The authors are grateful to Maja Kabiljo for her remarks and discussions that helped to improve the article. We are also sincerely thankful to an anonymous referee for the constructive remarks and references.

\end{document}